\documentclass[11pt,reqno]{amsart}
\usepackage{fourier}
\usepackage{fullpage}
\usepackage{mathrsfs,amssymb,graphicx,verbatim,amsmath,amsfonts}
\usepackage{paralist}
\usepackage[breaklinks,pdfstartview=FitH]{hyperref}
\usepackage{upgreek}
\usepackage[mathscr]{euscript}

\usepackage{xcolor}

\makeatletter
\def\resetMathstrut@{%
  \setbox\z@\hbox{%
    \mathchardef\@tempa\mathcode`\(\relax
    \def\@tempb##1"##2##3{\the\textfont"##3\char"}%
    \expandafter\@tempb\meaning\@tempa \relax
  }%
  \ht\Mathstrutbox@1.2\ht\z@ \dp\Mathstrutbox@1.2\dp\z@
}
\makeatother
\renewcommand{\le}{\leqslant}
\renewcommand{\ge}{\geqslant}
\renewcommand{\leq}{\leqslant}

\renewcommand{\setminus}{\smallsetminus}
\renewcommand{\gamma}{\upgamma}
\newcommand{\ud}[0]{\,\mathrm{d}}

\usepackage{dutchcal}

\newcommand{\GG}{\mathcal{G}}
\newcommand{\DD}{\mathscr{D}}
\newcommand{\sub}{\mathscr{C}}

\renewcommand{\det}{\mathrm{det}}

\newcommand{\sign}{\mathrm{sign}}

\newcommand{\n}{\{1,\ldots,n\}}

\newcommand{\X}{\mathbf X}
\newcommand{\Y}{\mathbf Y}

\newcommand{\f}{\varphi}

\renewcommand{\d}{\delta}

\newcommand{\e}{\varepsilon}
\newcommand{\R}{\mathbb R}
\newcommand{\Q}{\mathbb Q}

\newtheorem{theorem}{Theorem}

\newtheorem{remark}[theorem]{Remark}
\newcommand{\MM}{\mathcal{M}}

\newtheorem{conjecture}[theorem]{Conjecture}

\newtheorem{problem}[theorem]{Problem}

\newcommand{\sep}{\mathsf{SEP}}
\newcommand{\Vn}{V_{\!n}}
\newcommand{\SEP}{\mathsf{SEP}}

\renewcommand{\S}{\mathsf{S}}
\renewcommand{\subset}{\subseteq}
\newcommand{\C}{\mathbb C}

\DeclareMathOperator{\trace}{\mathrm{Trace}}

\newcommand{\J}{\mathfrak{J}}

\newcommand{\N}{\mathbb N}
\newcommand{\Z}{\mathbb Z}

\newcommand{\eqdef}{\stackrel{\mathrm{def}}{=}}

\DeclareMathOperator{\diam}{diam}

\newcommand{\A}{\mathsf{A}}
\renewcommand{\emptyset}{\varnothing}
\newcommand{\ee}{\mathsf{e}}
\newcommand{\conv}{\mathrm{conv}}

\newcommand{\iq}{\mathrm{iq}}
\newcommand{\Id}{\mathsf{I}}
\newcommand{\Idn}{\Id_{\!n}}
\newcommand{\M}{\mathsf{M}}

\newcommand{\SL}{\mathsf{SL}}
\renewcommand{\O}{\mathsf{O}}
\newcommand{\On}{\mathsf{O}_{\!n}}

\newcommand{\Part}{\mathscr{P}}
\newcommand{\vol}{\mathrm{vol}}

\newcommand{\HH}{\mathbb{H}}

\newcommand{\bE}{\mathbf{E}}

\newcommand{\prob}{\mathbf{Prob}}

\newcommand{\As}{\mathsf{A}^{\!\!*}}

\begin{document}

\title{The separation modulus of unitarily invariant matrix norms}

\author{Mustafa Alper Gunes}
\address{Mathematics Department\\ Princeton University\\ Fine Hall, Washington Road, Princeton, NJ 08544-1000, USA}
\email{magunes@princeton.edu}

\author{Assaf Naor}
\address{Mathematics Department\\ Princeton University\\ Fine Hall, Washington Road, Princeton, NJ 08544-1000, USA}
\email{naor@math.princeton.edu}

\thanks{A.~N.~was supported  by NSF grant DMS-2453936, BSF grant 2018223, and a Simons Investigator award. }

\maketitle

\vspace{-0.25in}

\begin{abstract} If $\X=(\M_n(\R),\|\cdot\|_\X)$ is a unitarily invariant normed space, namely, $\|\mathsf{U}\A\mathsf{V}\|_\X=\|\A\|_\X$ for every matrix $\A\in \M_n(\R)$ and every two  orthogonal matrices $\mathsf{U},\mathsf{V}\in \On$, then we prove that the spectral gap $\lambda(\X)$ of the Laplacian with Dirichlet boundary conditions on the unit ball $B_\X$ of $\X$ satisfies the asymptotic equivalence
\begin{equation}\label{eq:lambda in abstract}
\lambda(\X)\asymp n^3 \|\Id_n\|_\X^2, 
\end{equation}
where $\Idn$ is the $n$-by-$n$ identity matrix. This leads to a confirmation of the {weak isomorphic reverse isoperimetry conjecture} for $\X$, namely, we demonstrate that there exists a convex body $L=L_\X\subset B_\X$ such that 
$$
\big(\vol_{n^2}(L)\big)^{\frac{1}{n^2}}\asymp \big(\vol_{n^2}(B_{\X})\big)^{\frac{1}{n^2}},
$$
yet its isoperimetric quotient $\iq(L)=\vol_{n^2-1}(\partial L)/\vol_{n^2}(L)^{(n-1)/n}$ is at most a universal constant multiple of $n\asymp\iq(B_{\!\!\ell_{\!\!2}^{n^2}})$. We prove~\eqref{eq:lambda in abstract} through exact computations for a Jacobi orthogonal random matrix ensemble. 

As a corollary (and motivation)  of the above results, we deduce that the separation modulus of $\X$ satisfies
\begin{equation}\label{eq:in abstract}
\sep(\X)\asymp \sqrt{n}\|\Idn\|_\X\diam(B_\X),
\end{equation}
where  $\diam(B_\X)$ is the diameter of $B_\X$ with respect to the standard Euclidean  metric  on $\M_n(\R)$. Assuming oracle access to norm evaluations in $\X$, by combining~\eqref{eq:in abstract} with a new deterministic algorithm for efficiently computing a $O(1)$-approximation of the diameter of convex bodies in $\R^n$ that are given by a weak membership oracle and are symmetric with respect to coordinate permutations and reflections about the standard axes (this is  known to be impossible in the absence of such symmetries), we obtain an oracle polynomial time algorithm whose output is  guaranteed to be the separation modulus of  $\X$ up to positive universal constant factors. Another example of a consequence of~\eqref{eq:in abstract} is that for each $m\in \n$ the separation modulus of the $m$'th Ky Fan norm  on $\M_n(\R)$ is bounded from above and from below by positive universal constant multiples of $m\sqrt{n}$ if $m\ge \sqrt{n}$, and of $n$ if $m\le \sqrt{n}$. We also deduce from~\eqref{eq:in abstract} an upper bound on the Lipschitz extension modulus of $\X$ that improves over the previously best-known bound even in the special case when $\X$ is $\M_n(\R)$ equipped with the $\ell_{\!\!2}^n\to \ell_{\!\!2}^n$ operator norm. 
\end{abstract}

\section{Introduction}

Given $n\in \N$, a norm $\|\cdot\|$ on the space $\M_n(\R)$ of $n$-by-$n$ matrices with real entries is said to be unitarily invariant if $\|\mathsf{U}\A\mathsf{V}\|=\|\A\|$ for every $\A\in \M_n(\R)$ and every two matrices 
$\mathsf{U},\mathsf{V}$ that belong to the orthogonal group $\On$.\footnote{One defines analogously unitarily invariant norms on the spaces $\M_n(\C)$ and $\M_n(\mathbb{H})$ of $n$-by-$n$ matrices whose  entries belong to the complex numbers $\C$ or, respectively, the quaternions $\mathbb{H}$ (thus, the dimension over $\R$ of $\M_n(\C)$ equals $2n^2$ and the dimension over $\R$ of $\M_n(\mathbb{H})$ equals $4n^2$), but all of the results that are derived herein follow formally for complex and quaternionic entries from their counterparts for real entries; these simple deductions are carried out in Section~\ref{section:other entries} below. Alternatively, it is  possible to  perform  the  proofs  mutatis mutandis directly for complex and quaternionic matrices. For this reason, the considerations herein (excluding Section~\ref{section:other entries}) will treat only the case of matrices whose entries are in $\R$. } Norms with this property~\cite{vNe37,Sch50} form a rich class of ways to measure the size of matrices in an intrinsic manner that does not depend on the choice of orthogonal basses of $\R^n$ with respect to which the matrix representation of the corresponding linear operator is formed.  As such, unitarily invariant norms are of great importance and utility to a wide range of disciplines, ranging from multiple areas of pure mathematics and theoretical physics, to algorithm design, signal processing, and data analysis.  

Examples of unitarily invariant norms include the Schatten von-Neumann trace class $\S_p^n$ for every $1\le p\le \infty$, which are given by setting
$$
\forall \A\in \M_n(\R),\qquad \|A\|_{\S_p^n} \eqdef \Big(\trace \big((\A\A^{\!\!*})^{\frac{p}{2}}\big)\Big)^{\frac{1}{p}}=\Big(\sum_{i=1}^n s_i(\A)^p\Big)^{\frac{1}{p}}, 
$$
where $s_1(\A)\ge \ldots\ge s_n(\A)\ge 0$ are the singular values of $\A$, i.e., they are the  (decreasing rearrangement of the) eigenvalues  of the positive semidefinite matrix  $|\A|=(\A\As)^{1/2}$. If $p=\infty$, then  $\S_\infty^n$ coincides with the space of $n$-by-$n$ matrices equipped with their operator (a.k.a.~spectral) norm, when they are viewed as linear operators from $\ell_{\!\!2}^n$ to itself, i.e., 
$$\|\A\|_{\S_\infty^n}=\max_{x\in B_{\ell_{\!\!2}^n}} \|\A x\|_{\!\ell_{\!\!2}^n}=s_1(\A),$$ 
where  $B_\X=\{x\in \X:\ \|x\|_\X\le 1\}$ denotes the unit ball of a normed space $(\X,\|\cdot\|_\X)$. Often  $\S_1^n$ is called the nuclear norm on $\M_n(\R)$, and $\S_2^n$ is the standard Euclidean (a.k.a.~Hilbert--Schmidt or Frobenius) norm on $\M_n(\R)$, i.e., if $\A=(a_{ij})\in \M_n(\R)$, then $\|\A\|_{\S_2^n}=(\sum_{i=1}^n\sum_{j=1}^n a_{ij}^2)^{1/2}$. Additional  examples  are furnished~\cite{Fan51} by the Ky Fan norms $\|\cdot\|_{\mathsf{K}_1},\ldots,\|\cdot\|_{\mathsf{K}_n}$ that are defined as follows:
\begin{equation}\label{eq:def ky fan}
\forall m\in [n]\eqdef\n,\ \forall \A\in \M_n(\R),\qquad   \|A\|_{\mathsf{K}_m^n} \eqdef \sum_{i=1}^m s_i(\A).  
\end{equation}
More generally, let $\bE=(\R^n,\|\cdot\|_\bE)$ be a symmetric normed space, i.e., the following requirement holds:
\begin{equation}\label{eq:1 symmetric}
\forall x=(x_1,\ldots,x_n)\in \R^n,\ \forall \pi\in S_n, \ \forall (\e_1,\ldots,\e_n\in \{-1,1\},\qquad \|(\e_1x_{\pi(1)},\ldots,\e_n x_{\pi(n)})\|_\bE=\|x\|_\bE.
\end{equation}
The (real) unitary ideal $\S_\bE^n$ of $\bE$ is defined to be $\M_n(\R)$ equipped with the norm that is given by setting: 
\begin{equation}\label{eq:def ideal}
\forall \A\in \M_n(\R),\qquad \|\A\|_{\S_\bE^n}= \big\|\big(s_1(\A),\ldots,s_n(\A)\big)\big\|_\bE.
\end{equation}
A norm $\|\cdot\|$ on $\M_n(\R)$ is unitarily invariant if and only if there exists a symmetric norm $\|\cdot\|_\bE$ on $\R^n$ such that $\|\cdot\|=\|\cdot\|_{\S_\bE^n}$; see e.g.~\cite{Mir60} or~\cite[Chapter~IV]{Bha97} for a justification of this classical fact.

The purpose of the present article is to obtain optimal  randomized clustering of arbitrary unitarily invariant norms, which is a useful tool for a range of algorithmic tasks (see below), as well as for applications in mathematics (notably to metric embeddings~\cite{Bar98,FRT04,NR25} and Lipschitz extension~\cite{LN05,naor2022}).

Let $(\MM,d_\MM)$ be a metric space. For a partition $\Part$ of $\MM$ and a point $x\in \MM$, the cluster of $\Part$ to which $x$ belongs is denoted $\Part(x)$.  One says that $\Part$ is $\Delta$-bounded  for some $\Delta>0$ if $\diam_\MM(\Part(x))\le \Delta$ for every $x\in \MM$, where $\diam_\MM(\sub)=\sup\{d_\MM(y,z):\ y,z\in \MM\}$ denotes the $d_\MM$-diameter of $\emptyset\neq \sub\subset \MM$. Given $\sigma\ge 0$,  a random $\Delta$-bounded partition $\Part$  of $(\MM,d_\MM)$ is said to be $\sigma$-separating if:
\begin{equation}\label{eq:separation def}
          \forall x,y\in \MM,\qquad   \prob \big[\Part(x) \neq \Part(y)\big] \le  \frac{\sigma}{\Delta}d_\MM(x, y).
        \end{equation}
        Observe that if $\MM$ is a normed space, which is the only situation that will be investigated herein, then by dilation it suffices to consider only the special case $\Delta=1$ of~\eqref{eq:separation def}.  The separation modulus of $(\MM,d_\MM)$, which is denoted $\SEP(\MM,d_\MM)$ or simply $\SEP(\MM)$ when the underlying metric $d_\MM$ is clear from the context, is defined to be the infimum of those $\sigma\ge 0$ such that for every $\Delta>0$ there exists a distribution $\Part=\Part_{\!\!\Delta}$ over random partitions of $\MM$ that is both $\Delta$-bounded and $\sigma$-separating.

        \begin{remark}{\em 
Formally, when $(\MM,d_\MM)$ is an infinite metric space, one must impose some measurability requirement for the above definition to make sense: at the very least one has to demand that the event that appears in the left hand side of~\eqref{eq:separation def} is measurable with respect to the probability measure $\prob$ on the space of partitions of $\MM$. For the purpose of the present article, it suffices to assume the latter ``minimal'' measurability, namely, that the event $\{\Part(x)\neq \Part(y)\}$ is $\prob$-measurable for every $x,y\in \MM$. However, it is beneficial for applications to impose  stronger  measurability assumptions; such a foundational treatment of random partitions of infinite spaces has been carried out in~\cite{LN05,naor2022}. One of the applications of our main result (Theorem~\ref{thm:sep formula}  below) is to deduce an improved Lipschitz extension result (Theorem~\ref{thm:ext cor} below), for which the aforementioned stronger  measurability assumptions are needed. However, this will be achieved by substituting Theorem~\ref{thm:sep formula} into the framework of~\cite{naor2022}, which is shown in~\cite{naor2022} to indeed yield the desired measurability. Thus, we do not need to recall the measurability requirements  as they do not directly pertain to the ensuing discussion. Alternatively, our results are meaningful and new even if one only treats random partitions of finite subsets of unitarily invariant matrix norms rather than random partitions of $\M_n(\R)$, a setting to which measurability considerations are irrelevant.  }
\end{remark}

The above notion of randomized clustering is due to~\cite{bartal1996probabilistic} (precursors that considered it implicitly, with a variety of applications, appear in~\cite{LR88,AP90,AKPW91,KPR93}); see e.g.~\cite{Nao17,naor2022} for more on the history. By~\cite{bartal1996probabilistic}  every $n$-point metric space $\MM$ satisfies $\SEP(\MM)=O(\log n)$, and there are metric spaces for which this bound is sharp. In~\cite{KR11},  a randomized polynomial time algorithm was designed that takes as input a finite metric space $\MM$ and outputs a factor $2$ approximation to $\sep(\MM)$. See~\cite{Bar98,ccg98,KRSS01,FRT04,LN05, KR11,BFKMNNS14,Nao17,naor2022,krauthgamer2025lipschitz,NR25,KPS25} for an indication of the substantial literature on this topic. The investigation of separating random partitions of finite-dimensional normed spaces, which is what we study herein, originates in~\cite{PR98},   motivated by applications to network routing and distributed computing. The  work~\cite{ccg98} sharpened and generalized the bounds of~\cite{PR98}, and influenced later research, e.g.~\cite{LN05,AI06,naor2022}; similar partitioning schemes appeared implicitly in the earlier work~\cite{KMS98} on SDP-based algorithms for graph colorings.

While the separation modulus  is  an intricate nonlinear invariant that could be difficult to evaluate because  its definition consists of an optimization over  all probability distributions on arbitrary $\Delta$-bounded partitions into clusters that may be very complicated,  here we will obtain a simple asymptotic formula (up to universal constant factors) for  the separation modulus of any unitarily invariant matrix norm:

\begin{theorem}\label{thm:sep formula} Let $\X=(\M_n(\R),\|\cdot\|_\X)$ be a unitarily invariant normed space. Then\footnote{In addition to  the usual $O(\cdot),o(\cdot),\Omega(\cdot), \Theta(\cdot)$ notation, we will use throughout this text  the following conventions for asymptotic notation: Given $a,b>0$, by writing
$a\lesssim b$ or $b\gtrsim a$ we mean that $a\le \kappa b$ for some
universal constant $\kappa>0$, and $a\asymp b$
stands for $(a\lesssim b) \wedge  (b\lesssim a)$; when we will need to allow for dependence on parameters, we will indicate it by subscripts.}    
\begin{equation}\label{eq:in thm}
\sep(\X)\asymp \sqrt{n}\|\Idn\|_\X\diam(B_\X),
\end{equation}
where  $\diam(B_\X)= \diam_{\S_2^n}(B_\X)$ is the diameter with respect to the standard Euclidean metric on $\M_n(\R)$  of the unit ball $B_\X$ of $\X$, and $\Id_n\in \M_n(\R)$ is the $n$-by-$n$ identity matrix. 
\end{theorem}

The novel content of Theorem~\ref{thm:sep formula} is the upper bound on $\SEP(\X)$ in~\eqref{eq:in thm}, i.e., Theorem~\ref{thm:sep formula}  provides an improved randomized clustering of $\X$. The matching impossibility result is due to~\cite[Corollary~79]{naor2022}, where it was proved that $\sqrt{n}\|\Idn\|_\X\diam(B_\X)\lesssim \sep(\X)\lesssim \sqrt{n\log n}\|\Idn\|_\X\diam(B_\X)$ in the setting of Theorem~\ref{thm:sep formula}. Therefore, what Theorem~\ref{thm:sep formula} achieves  is removing the last   redundant lower-order term, thus obtaining an asymptotic evaluation of $\SEP(\X)$ that is optimal up universal constant factors. 

For the operator norm, the aforementioned bounds from~\cite{naor2022} are $n\lesssim \SEP(\S_\infty^n)\lesssim n\sqrt{\log n}$, and these were the previously best-known estimates on $\SEP(\S_\infty^n)$. By Theorem~\ref{thm:sep formula},   we now know that $\SEP(\S_\infty^n)\asymp n$. 

As another example, it is straightforward to check from~\eqref{eq:def ky fan} that for each $m\in [n]$ the $\S_2^n$-diameter  of the unit ball of the $m$'th Ky Fan norm $\mathsf{K}_m^n$ equals $2\sqrt{n+m(m-1)}/m\asymp \max\{1,\sqrt{n}/m\}$, so Theorem~\ref{thm:sep formula}  demonstrates that $\SEP(\mathsf{K}_m^n)$ is bounded from above and from below by positive universal constant multiples of $m\sqrt{n}$ if $m\ge \sqrt{n}$, and of $n$ if $m\le \sqrt{n}$, thus providing an interpolation between the operator  norm $\mathsf{K}_1^n$  and the nuclear norm for $\mathsf{K}_n^n$ that somewhat curiously  stabilizes (up to universal constant factors) at $m\asymp \sqrt{n}$.

Theorem~\ref{thm:sep formula} yields progress on the important and extensively-studied classical topic of extending Lipschitz functions; see e.g.~\cite{BB12,naor2022} and the references therein. The Lipschitz extension modulus of a metric space $(\MM,d_\MM)$, denoted $\ee(\MM,d_\MM)$ or simply $\ee(\MM)$ when the underlying metric $d_\MM$ is clear from the context,  is the infimum over those $L\in [1,\infty]$ such that for every subset $\sub\subset \MM$, every Banach space $(\Y,\|\cdot\|_\Y)$, and every $1$-Lipschitz function $f:\sub\to \Y$, there exists an $L$-Lipschitz function $F:\MM\to \Y$ whose restriction to $\sub$ coincides with $f$. By~\cite{LN05} we have $\ee(\MM)\lesssim\SEP(\MM)$, and the forthcoming work~\cite{BN25} demonstrates that  $\ee(\MM)\gtrsim \sqrt{\SEP(\MM)}$, so $\ee(\MM)<\infty$ if and only if $\SEP(\MM)<\infty$. A substitution of the former upper bound on $\ee(\MM)$ provides the following extension theorem:

\begin{theorem}\label{thm:ext cor} If $\X=(\M_n(\R),\|\cdot\|_\X)$ is a unitarily invariant normed space, then
$
\ee(\X)\lesssim \sqrt{n}\|\Idn\|_\X\diam(B_\X).
$
\end{theorem}

The special case $\X=\S_\infty^n$ of Theorem~\ref{thm:ext cor} asserts that for every subset $\sub$ of $\M_n(\R)$ and every Banach space $(\Y,\|\cdot\|_\Y)$, if  $f:\sub\to  \Y$ is $1$-Lipschitz with respect to the operator norm (and the given norm $\|\cdot\|_\Y$ on the target $\Y$), then it can be extended to a $\Y$-valued function that is defined on all of $\M_n(\R)$ and is $L$-Lipschitz with respect to the operator norm, where $L\lesssim n$. The previously best-known~\cite{naor2022} upper bound on $L$ here was $L\lesssim n\sqrt{\log n}$, and the corresponding best impossibility result that is currently available~\cite{naor2022} is that there exist suitably chosen $\sub,\Y,f$ as above for which necessarily $L\gtrsim \sqrt{n}$. Thus, we now know that: \begin{equation}\label{ext bounds}
{\textstyle\sqrt[4]{\dim(\S_\infty^n)}}=\sqrt{n}\lesssim \ee(\S_\infty^n)\lesssim n={\textstyle\sqrt{\dim(\S_\infty^n)}}.
\end{equation}

\begin{problem}{\em 
It remains an  interesting open question to determine the growth rate of $\ee(\S_\infty^n)$ as $n\to \infty$. For that matter, there is no unitarily invariant normed space $\X=(\M_n(\R),\|\cdot\|_\X)$ for which the value of $\ee(\X)$ is  known up to $O(1)$ factors, or even up to lower-order factors. Since $\diam(B_\X)\ge 2\|\Idn\|_{\S_2^n}/\|\Idn\|_\X=2\sqrt{n}\|\Idn\|_\X$, Theorem~\ref{thm:ext cor} does not provide an upper bound on $\ee_n(\X)$ that is $o(n)$. At the same time, the general lower bound on $\ee(\X)$ in~\cite[Theorem~96]{naor2022} is at most $\sqrt{n}$, since the Banach--Mazur distance between $\X$ and an $n^2$-dimesnional Hilbert space is at most $\sqrt{n}$, as seen by applying John's theorem~\cite{Joh48} to the symmetric space $\bE=(\R^n,\|\cdot\|_\bE)$ for which $\X=\S_\bE^n$. Thus, the bounds in~\eqref{ext bounds} are the best one could hope for using the current state of the art even when $\S_\infty^n$ is replaced by any unitarily invariant norm on $\M_n(\R)$. It seems that a substantial new idea will be needed in order to improve either the upper or lower estimate in~\eqref{ext bounds}. }\end{problem}

As another application of Theorem~\ref{thm:sep formula}, if we are given oracle access to (approximate) norm computations in a unitarily invariant normed space $\X=(\M_n(\R),\|\cdot\|_\X)$, then we obtain a (deterministic) algorithm that evaluates its separation modulus $\SEP(\X)$ up to universal constant factors in oracle polynomial time.\footnote{A weak membership oracle suffices here (one can obtain from such a membership oracle a norm evaluation oracle by binary search of possible upper bounds on the norm of the queried vector). See~\cite{GLS93} for the relevant (standard)  background in algorithmic convex geometry. In particular, ``oracle polynomial time'' means that  the algorithm consists of a binary Turing machine that is augmented by the given oracle, and its running time is polynomial in the usual sense, with the understanding that the time required by each call to the oracle is only what is needed to write the question onto, and read the answer from, a tape of the Turing machine. Thus, here we demand that the algorithm outputs a number that is guaranteed to be bounded from above and from below by positive universal constant multiples of $\SEP(\X)$, yet it performs only $n^{O(1)}$ norm evaluations.} 

Theorem~\ref{thm:sep formula} reduces the above task to devising an algorithm that outputs a $O(1)$-factor approximation to $\diam(B_\X)$ in oracle polynomial time.  Write $\X=\S_\bE^n$ for a symmetric normed space $\bE=(\R^n,\|\cdot\|_\bE)$. The $\S_2^n$-diameter of $B_\X$ is equal to the $\ell_{\!\!2}^n$-diameter of $B_\bE$. And, $\X$ has the desired oracle if and only if $\bE$ does.~\footnote{We will use here only the trivial direction of this equivalence, which is simply restricting the oracle for $\X$ to diagonal matrices. The  reverse direction incurs the $n^{O(1)}$- time overhead of computing the singular values of the query matrix.}   So, the goal is obtaining  an oracle polynomial time $O(1)$-factor approximation algorithm to $\diam_{\ell_{\!\!2}^n}(B_\bE)$.

 A  reduction to efficient Euclidean diameter estimation should give one pause, as  the latter task cannot be achieved for general convex bodies, i.e., there is no oracle polynomial time algorithm    that  outputs a $O(1)$-factor approximation to the $\ell_2^n$-diameter of an arbitrary convex body $K$ in $\R^n$ that is given by a weak membership oracle; for deterministic algorithms, this impossibility result was proved in~\cite{BF87}, and randomized algorithms were ruled out in~\cite{BGKKLM98,BGKKLS01}. This hardness of approximation is also proved in~\cite{BGKKLS01} for the $\ell_{\!\!p}^n$-diameter for every $1\le p<\infty$ (see Remark~\ref{rem:approx ratios}  below for a description of the inapproximability factors). Nevertheless, in the presence of additional symmetries, we will prove the following theorem, which yields the aforementioned algorithm for approximating $\SEP(\X)$ using  Theorem~\ref{thm:sep formula}:

\begin{theorem}\label{thm:alg rad} For every $n\in \N$, $D\ge 1$ and $0<\d\le \frac12$, there exists $d=d(D,\d)\in \N$ satisfying
\begin{equation}\label{eq:deg}
d\lesssim D+\frac{\log\left(\frac{1}{\d}\right)}{\d},
\end{equation}
such that  for every $1\le p\le \infty$  there exists a deterministic  algorithm $\mathrm{Alg}_p$ that takes as input $n,D,\e$ and a symmetric normed space $\bE=(\R^n,\|\cdot\|_\bE)$  whose unit ball is given by a weak membership oracle and satisfies
\begin{equation}\label{eq:norms are reasonable}
\frac{1}{e^{n^{D}}}\le \|e_1\|_\bE=\ldots=\|e_n\|_\bE\le e^{n^{D}}.
\end{equation}
The algorithm $\mathrm{Alg}_p$ outputs in oracle time $n^d$ a number $\mathrm{Alg}_p(n,D,\d,\bE)$ that satisfies 
\begin{equation}\label{eq:guarantees of alg}
(1-\d)\diam_{\ell_{\!\!p}^n}(B_\bE)\le \mathrm{Alg}_p(n,D,\d,\bE)\le (1+\d)\diam_{\ell_{\!\!p}^n}(B_\bE). 
\end{equation}
\end{theorem}

\subsection{Reverse Faber--Krahn}\label{sec:FB} The  new ingredient  of our proof of Theorem~\ref{eq:in thm} is   an asymptotic evaluation of the spectral gap of the Laplacian with Dirichlet boundary conditions on the unit ball $B_\X\subset \M_n(\R)$ of any unitarily invariant norm $\X=(\M_n(\R),\|\cdot\|_\X)$. We will next describe this;  all of the relevant (classical and rudimentary) background on the spectral properties of the Dirichlet Laplacian can be found in~\cite{PS51,Cha84}. 

Fix $n\in \{2,3,\ldots\}$.  For a bounded domain $\Omega\subset \R^n$,  let $\lambda(\Omega)$ be the smallest $\lambda>0$ such that there exists a function $f:\Omega\to \R$ that is smooth on the interior of $\Omega$, vanishes on the boundary $\partial \Omega$ of $\Omega$, and satisfies $\Delta f= - \lambda f$ on the interior of $\Omega$, where $\Delta=\sum_{i=1}^n\partial^2/(\partial x_i^2)$ is the standard Laplacian on $\R^n$. Set $\lambda(\Y)=\lambda(B_\Y)$ for each normed space $\Y=(\R^n,\|\cdot\|_\Y)$. Here, we obtain the following result: 

 \begin{theorem}\label{thm:spectral} For every $n\in \N$, every unitarily invariant normed space $\X=(\M_n(\R),\|\cdot\|_\X)$ satisfies:
\begin{equation}\label{eq:RFB in thereom}
\lambda(\X)\asymp n^3\|\Id_n\|_\X^2. 
\end{equation}
\end{theorem}

In the setting of Theorem~\ref{thm:spectral}, the  following asymptotic equivalence is a different way to express~\eqref{eq:RFB in thereom}:  
\begin{equation}\label{eq:invariant RFB}
\lambda(B_\X)\vol_{\dim(\X)}^{\frac{2}{\dim(\X)}}(B_\X)\asymp\dim(\X), 
\end{equation}
where  for $m\in \N$ and $0\le \alpha\le m$ the $\alpha$-dimensional Hausdorff measure on $\R^m$  is denoted  $\vol_\alpha (\cdot)$; the only values of $\alpha$ that occur herein are $\alpha\in \{m,m-1\}$  (Lebesgue measure and surface area measure, respectively).  To see that~\eqref{eq:invariant RFB} and~\eqref{eq:RFB in thereom} coincide, note that $\dim(\X)=n^2$ and  by~\cite[equation~(2.2)]{Sch82} we have: 
\begin{equation}\label{eq:volume of Sinfty}
\vol_{n^2} (B_\X)^{\frac{1}{n^2}}\asymp \frac{1}{\sqrt{n}\|\Id_n\|_\X}.
\end{equation} 

 The Faber--Krahn inequality~\cite{Fab23,Kra26} (originally conjectured in~\cite{Ray45}) states that the left hand side of~\eqref{eq:invariant RFB}  is minimized when $B_\X$ is replaced by the unit ball of the standard Euclidean metric (in our case $\S_2^n$, though the analogous statement holds if $\M_n(\R)$ is replaced by  $\R^m$ of any dimension).  By a direct computation  (see e.g.~\cite[page~32]{naor2022}), this  classical  extremal statement implies that the left hand side of~\eqref{eq:invariant RFB} is bounded from above and from below by positive universal constant multiples of $\dim(\X)$. 
 
 The new content of~\eqref{eq:invariant RFB} is therefore the reverse inequality, i.e., the estimate $\lambda(\X)\lesssim n^3\|\Id_n\|_\X^2$    in Theorem~\ref{thm:spectral}. This confirms for (unit balls of) unitarily invariant matrix norms the following much more general reverse Faber--Krahn phenomenon that was conjectured in~\cite{naor2022}:
\begin{conjecture}[reverse Faber--Krahn]\label{con:reverse FK} For every origin-symmertric convex body $K\subset \R^m$ (thus, $x\in K$ if and only if $-x\in K$) there is a volume-preserving linear transformation $S\in \SL_m(\R)$ with $\lambda(SK)\vol_{m}(K)^{2/m}\asymp m$. 
\end{conjecture}

 Theorem~\ref{thm:spectral} demonstrates that Conjecture~\ref{con:reverse FK} holds for any unitarily invariant norm on $\M_n(\R)\cong \R^{n^2}$ with the matrix $S$ being the ($n^2$-by-$n^2$) identity on $\M_n(\R)$; this too was conjectured in~\cite{naor2022} (see specifically Conjecture~50 there, as well as~\cite[Problem~44]{naor2022} and~\cite[Remark~172]{naor2022}).

\subsection{Weak isomorphic reverse isoperimetry}\label{sec:weak reverse} It was proved in~\cite{naor2022} that Theorem~\ref{thm:spectral} implies Theorem~\ref{thm:sep formula}. The link between these two results is Theorem~\ref{thm:weak revers iso unitarily} below. For its formulation, recall that the isoperimetric quotient of a convex body $K\subset \R^m$  is defined to be $\iq(K)=\vol_m(K)/\vol_{m-1}(\partial K)^{1-1/m}$.

\begin{theorem}\label{thm:weak revers iso unitarily} For every $n\in \N$, if  $\X=(\M_n(\R),\|\cdot\|_\X)$ is a unitarily invariant normed space, then there exists an origin-symmetric convex body $L\subset B_\X$ such that: 
\begin{equation}\label{eq:weak revese iso for unitarily}
\vol_{n^2}(L)^{\frac{1}{n^2}}\gtrsim \vol_{n^2}(B_\X)^{\frac{1}{n^2}}\qquad\mathrm{yet}\qquad  \iq(L)\lesssim n. 
\end{equation} 
\end{theorem}

It was proved  in~\cite[Section~1.6.2]{naor2022} that for each unitarily invariant normed space $\X=(\M_n(\R),\|\cdot\|_\X)$, the validity of Theorem~\ref{thm:spectral}  for $\X$ is equivalent to that validity of Theorem~\ref{thm:weak revers iso unitarily} for $\X$; see specifically the explanation why~\cite[Conjecture~49]{naor2022} and~\cite[Conjecture~50]{naor2022} are equivalent that appears in~\cite[page~43]{naor2022}, for which equation~(1.62) of~\cite{naor2022} is  crucial input, as well as the convexity and uniqueness of the so-called Cheeger body (see page~31 of~\cite{naor2022}), per~\cite{AC09}. Therefore, by proving Theorem~\ref{thm:spectral}  we will also establish Theorem~\ref{thm:weak revers iso unitarily} , but we will next briefly explain the meaning of Theorem~\ref{thm:weak revers iso unitarily}  and its utility for proving Theorem~\ref{thm:sep formula}.

For $m\in \N$, the isoperimetric theorem implies that  the isoperimetric quotient of a convex body $K\subset \R^m$ is at least the isoperimetric quotient of the $\ell_2^m$-ball. Direct computation shows that the latter is   bounded from above and from below by positive universal constant multiples of $\sqrt{m}$. Thus, Theorem~\ref{thm:weak revers iso unitarily} says that one can inscribe in $B_\X$ an origin-symemtric convex body $L$ that is large in the sense that the first inequality in~\eqref{eq:weak revese iso for unitarily} holds, yet its isoperimetric quotient is as small as possible up to universal constant factors, i.e., $\iq(L)$ has the same  order of magnitude as the isoperimetric quotient of the Hilbert--Schmidt unit ball in $\M_n(\R)$.   This confirms for (unit balls of) unitarily invariant matrix norms the following much more general  phenomenon that was conjectured in~\cite{naor2022}:

\begin{conjecture}[weak isomorphic reverse isoperimetry~\footnote{The adjective ``weak'' was used in~\cite{naor2022} to describe this conjecture  because  in~\cite{naor2022} a stronger phenomenon was also conjectured; we do not need to recall that stronger conjecture because  the present work does not address it, and furthermore its weaker variant that we establish herein suffices for proving Theorem~\ref{thm:sep formula}.}]\label{conj:weak reverse iso} For every origin-symmetric convex body $K\subset \R^m$ there exist $S\in \SL_m(\R)$ and an origin-symemtric convex body $L\subset SK$ that satisfies: 
$$
\vol_{m}(L)^{\frac{1}{m}}\gtrsim \vol_{m}(B_\X)^{\frac{1}{m}}\qquad\mathrm{yet}\qquad  \iq(L)\lesssim \sqrt{m}. 
$$
\end{conjecture}

Analogously to the discussion in Section~\ref{sec:FB}, Theorem~\ref{thm:weak revers iso unitarily} demonstrates that Conjecture~\ref{conj:weak reverse iso} holds for any unitarily invariant norm on $\M_n(\R)$ with the matrix $S$ being the identity on $\M_n(\R)$; this too was conjectured in~\cite{naor2022} (see specifically Conjecture~49 there).

By~\cite{naor2022}, the upper bound on $\SEP(\X)$ of Theorem~\ref{thm:sep formula} (which, as we explained earlier, is the part of Theorem~\ref{thm:sep formula} that remained to be proved because the matching lower bound is due to~\cite{naor2022}) follows from Theorem~\ref{thm:weak revers iso unitarily}; specifically, combine~\cite[Corollary~79]{naor2022} with equation~(2.2) of~\cite{Sch82}, which is reproduced as equation~(6.131) in~\cite{naor2022}. The resulting construction in~\cite{naor2022}  of the random partition of $\X$ can be illustrated informally as follows. Let  $L\subset B_\X$ be the auxiliary convex body that  Theorem~\ref{thm:sep formula} provides. Iteratively remove random\footnote{One must  define a suitable distribution of those random shifts here; for the purpose of the present sketch, one could aim instead to obtain a separated partition of an arbitrarily large bounded domain in $\M_n(\R)$, in which case the shifts can be chosen Lebesgue-uniformly at random.} translates of $L$ from $\M_n(\R)$; this yields almost surely a random partition of $\M_n(\R)$ into sets (translates of $L$ minus the union of those translates from the preceding steps of this random iterative removal)  that are of $\X$-diameter at most $1$, as $L\subset B_\X$, and~\cite{naor2022} demonstrates that the desired estimate on the separation probability holds thanks to the upper bound on the isoperimetric quotient of $L$ in~\eqref{eq:weak revese iso for unitarily}. 

\begin{remark}
{\em While the above procedure yields a random partition of $\X$ whose separation modulus is optimal up to universal constant factors, it is not entirely constructive. From the algorithmic perspective it suffices to bound $\SEP(\X)$ from above, per Theorem~\ref{thm:alg rad} and the discussion that precedes its statement, as by~\cite{KR11} there is a randomized polynomial time algorithm that takes as input any finite subset $\sub$ of $\X$ and outputs a $2\sep(\X)$-separating partition of $\sub$. One could hope to extract more algorithmic information by assuming that $B_\X$ has a weak membership oracle and aiming to obtain  a convex body $L$ as in  Theorem~\ref{thm:sep formula} which also has a weak membership oracle that is allowed to make $n^{O(1)}$ calls to the oracle of $B_\X$ and $n^{O(1)}$ additional Turing machine operations. We do not currently have this because the body $L$  of Theorem~\ref{thm:sep formula} is defined implicitly as the solution of a variational problem. We expect that the aforementioned  task could be  achieved by computing in oracle polynomial time  a sufficiently good approximation to the first nontrivial eigenfunction $\f:B_\X\to \R$ of the Dirichlet Laplacian on $B_\X$, i.e., $\f$ is smooth on the interior of $B_\X$, it vanishes on $\partial B_\X$, and $-\Delta\f=\lambda(\X)\f$ point-wise. By~\cite{BL76} we know that $\f$ is a log-concave function, whence its sub-level sets are convex. We could then repeat the above procedure with $L$ replaced a suitable sub-level set of $\f$. We conjecture that this will yield an optimal separating partition of  (a large bounded domain of) $\X$, and it remains open to implement the above strategy efficiently. 

The above proposal is an (especially interesting) instance of a general program to estimate efficiently (oracle polynomial time) solutions of geometrically meaningful partial differential equations; another important example is the first nontrivial eigenfunction of the  Laplacian on a convex body with Neumann boundary conditions, which controls the Kannan--Lov\'asz--Simonovits conjecture~\cite{KLS95}; see~\cite{Mil09,LV19}). }
\end{remark}

\section{Proof of Theorem~\ref{thm:spectral}  }

Fix $n\in \N$. It was proved in~\cite[Remark~172]{naor2022} that if Theorem~\eqref{thm:weak revers iso unitarily} holds for $\S_\infty^n$, then Theorem~\eqref{conj:weak reverse iso} holds for every unitarily invariant normed space $\X=(\M_n(\R),\|\cdot\|_\X)$. By the equivalence of Theorem~\ref{thm:spectral}  and Theorem~\ref{thm:weak revers iso unitarily} that was proved in~\cite{naor2022} and we discussed in Section~\ref{sec:weak reverse}, it therefore suffices to prove Theorem~\ref{thm:spectral} when $\X=\S_\infty^n$, where we recall that $\S_\infty^n$ is  $\M_n(\R)$ equipped with the $\ell_{\!\!2}^n\to \ell_{\!\!2}^n$ operator norm. So, our goal in the rest of this section will be to demonstrate that $\lambda(\S_\infty^n)\asymp n^3$. Only the upper bound 
\begin{equation}\label{eq:goal operator norm}
\lambda(\S_\infty^n)\lesssim n^3
\end{equation}
remains to be proved, thanks to the Faber--Krahn inequality,  as we explained in Section~\ref{sec:FB}.

For every $p,n\in \N$ define a function  $f_n^p:B_{\S_\infty^n}\to [0,\infty)$ by setting:
\begin{equation}
 \forall \A\in B_{\S_\infty^n},\qquad    f_n^p(\A)\eqdef \det\big(\Id_n-\A \As\big)^p=  \prod_{k=1}^n \big(1-s_k(\A)^2\big)^{p}.
    \label{testfxndef}
\end{equation}
Then, $f_n^p$ is smooth (in fact, it is a polynomial in the entries of $\A$), and the second equality in~\eqref{testfxndef} makes it evident that $f_n^p$ vanishes on $\partial B_{\S_\infty^n}=\{\A\in \M_n(\R):\ s_1(\A)=1\}$. The standard (see e.g.~\cite{Cha84}) Rayleigh quotient characterization of the smallest nonzero eigenvalue of $(-\Delta)$ therefore gives: 
\begin{equation}\label{eq:bound by quotient}
\lambda(\S_\infty^n)\le \frac{\int_{B_{\S_\infty^n}}\|\nabla f_n^p(\A)\|_{\S_2^n}^2\ud \A}{\int_{B_{\S_\infty^n}} |f_n^p(\A)|^2\ud \A}.
\end{equation}

Theorem~\ref{thm:spectral}  will thus be proven if we will evaluate exactly  the right hand side of~\eqref{eq:bound by quotient} as follows:
\begin{equation}\label{eq:ratio exactly}
\frac{\int_{B_{\S_\infty^n}}\|\nabla f_n^p(\A)\|_{\S_2^n}^2\ud \A}{\int_{B_{\S_\infty^n}} f_n^p(\A)^2\ud \A}=\frac{pn^2 (4 p+n)^2 }{2 (2 p-1) (4 p+1)}.
\end{equation}
Indeed, the right hand side of~\eqref{eq:ratio exactly} is bounded from above and from below by positive universal constant multiples of $n^2\max\{n^2/p,p\}$. Therefore, up to universal constant factors  setting $p=n$ minimizes the right hand side of~\eqref{eq:ratio exactly}. For this choice of $p$, the Rayleigh quotient estimate~\eqref{eq:bound by quotient} becomes $\lambda(\S_\infty^n)\lesssim n^3$, which is our current goal~\eqref{eq:goal operator norm}. Observe that the above reasoning shows that any choice of $p$ other than $p\asymp n$ does not yield the desired estimate~\eqref{eq:goal operator norm} .  The rest of this section will be devoted to establishing~\eqref{eq:ratio exactly}.

Given $n\in \N$, we will denote by $\Vn:\R^n\to [0,\infty)$ the following (Vendermunde determinant) function:
\begin{equation}\label{eq:def V}
\forall x=(x_1,\ldots,x_n)\in \R^n,\qquad \Vn(x)\eqdef \prod_{i=1}^{n-1}\prod_{j=i+1}^n|x_i-x_j|.
\end{equation} 
Define a probability density $\J_n^p:[0,1]^n\to [0,\infty)$ by setting:
\begin{equation}\label{eq:J density def}
\forall x=(x_1,\ldots,x_n)\in [0,1]^n,\qquad \J_n^p(x) \eqdef \frac{\Vn(x)}{Z_n(p)}\prod_{r=1}^n\frac{x_r^{2p}}{\sqrt{1-x_r}},
\end{equation}
where $Z_n(p)>0$ denotes the normalization factor such that $\int_{[0,1]^n}\J_n^p(x)\ud x=1$, i.e., 
\begin{equation}\label{normalization formula}
Z_n(p)\eqdef \int_{[0,1]^n} \Vn(x)\prod_{r=1}^n\frac{x_r^{2p}}{\sqrt{1-x_r}}\ud x=\frac{n!}{\pi^{\frac{n}{2}}}\prod_{s=1}^n\frac{\Gamma\big(2m+\frac12(s+1)\big) \Gamma\big(\frac12s\big)^2}{\Gamma\big(2m+\frac12(n+s+1)\big) }.
\end{equation}
The second equality in~\eqref{normalization formula} is a special case of the classical Selberg integral formula~\cite{Sel44}. The parameter $Z_n(p)$ is commonly denoted in the literature (see e.g.~\cite{aomoto}) as $S_n(\alpha,\beta,\gamma)$ for the setting of parameters $(\alpha,\beta,\gamma)=(2p+1,1/2,1/2)$; using the ad hoc shorter notation $Z_n(p)$ is beneficial here since we will not need to consider probability  densities other than $\J_n^p$. Also, we chose to use the letter $\J$ for the probability density on $[0,1]^n$  in~\eqref{eq:J density def} to indicate that it belongs to the $3$-parameter family of densities of Jacobi random $n$-by-$n$  matrix ensembles; specifically, using the parametrization of~\cite[equation~(4.1.5)]{anderson_guionnet_zeitouni_2009}, the relevant setting of parameters in the present case are $(\beta,r,s)=(1,4p+1,0)$.

There exists $c_n>0$ such that for every  Borel-measurable function $f:[0,\infty)^n\to [0,\infty)$ that is invariant under permutations of the coordinates, i.e., $f(s_{\pi(1)},\ldots,s_{\pi(n)})=f(s)$ for every  permutation $\pi\in S_n$ and every point $s=(s_1,\ldots,s_n)\in [0,\infty)^n$, the following identity holds: 
\begin{equation}\label{eq:weyl}
        \int_{\M_n(\R)} f\big(s_1(\A),\ldots,s_n(\A)\big) \ud A = c_n \int_{[0,\infty)^n}  f(s) \prod_{i=1}^{n-1}\prod_{j=s+1}^n \big|s_i^2 -s_j^2\big| \ud s.
    \end{equation}
\eqref{eq:weyl}  is a  special case of the Weyl integration formula~\cite{Wey39}; see e.g.~\cite[Proposition~4.1.3]{anderson_guionnet_zeitouni_2009}. The exact value of the parameter $c_n$ will not have any role in the ensuing reasoning (it will cancel out in our computations), but it is worthwhile to note in passing that $c_n^{1/n^2}\asymp 1/\sqrt{n}$; see e.g.~\cite{For10} for much more on such matters.  

Thanks to~\eqref{eq:weyl}, the denominator in the left hand side of~\eqref{eq:ratio exactly} can be evaluated as follows:
\begin{equation}\label{eq:first change of variable}
\int_{B_{\S_\infty^n}} f_n^p(\A)^2\ud \A\stackrel{\eqref{testfxndef}\wedge \eqref{eq:weyl}}{=} c_n\int_{[0,1]^n}\prod_{r=1}^n \big(1-s_r^2\big)^{2p}\prod_{i=1}^{n-1}\prod_{j=s+1}^n \big|s_i^2 -s_j^2\big| \ud x\stackrel{\eqref{normalization formula}}{=}\frac{c_n}{2^n}Z_n(p), 
\end{equation}
where for the last step of~\eqref{eq:first change of variable} we performed the (diagonal) change of variable $s=(\sqrt{1-x_1},\ldots,\sqrt{1-x_n})$, whose Jacobian is $(-1)^n/(2^n\prod_{r=1}^n\sqrt{1-x_r})$.

Next, we need to understand the Hilbert--Schmidt norm of the gradient of  $(\A\in B_{\S_\infty^n})\mapsto f_n^p(\A)$. For this, as a special case of~\cite[Theorem~7.1]{Lewis2005NonsmoothAO} we see that if we define $g_n^p:\R^n\to \R$ by 
\begin{equation}\label{eq:def our diag g}
\forall x\in \R^n,\qquad g_n^p(x)\eqdef \prod_{i=1}^n \big(1-x_i^2\big)^{p},
\end{equation}
then for every matrix $\A$ in the interior of $B_{\S_\infty^n}$ there exist $\mathsf{U}=\mathsf{U}_{\A,n,p},\mathsf{V}=\mathsf{V}_{\A,n,p}\in \O_n$ such that
\begin{equation*}
\nabla f_n^p(\A)=\mathsf{U} \begin{pmatrix} \frac{\partial g_n^p}{\partial x_1}\big(s_1(\A),\ldots,s_n(\A)\big) & 0 &  \dots&0 \\
  0 & \frac{\partial g_n^p}{\partial x_2}\big(s_1(\A),\ldots,s_n(\A)\big)&  \ddots & \vdots\\
            \vdots & \ddots & \ddots &0\\
              0 & \dots & 0&\frac{\partial g_n^p}{\partial x_n}\big(s_1(\A),\ldots,s_n(\A)\big)
                       \end{pmatrix}\mathsf{V}.
\end{equation*}
 By substituting the definition~\eqref{eq:def our diag g} of $g_n^p$, we thus have:
 $$
\nabla f_n^p(\A) = 2p\prod_{i=1}^n \big(1-s_i(\A)^2\big)^{p}\mathsf{U} \begin{pmatrix} \frac{s_1(\A)}{1-s_1(\A)^2}& 0 &  \dots&0 \\
  0 & \frac{s_2(\A)}{1-s_2(\A)^2}&  \ddots & \vdots\\
            \vdots & \ddots & \ddots &0\\
              0 & \dots & 0&\frac{s_n(\A)}{1-s_n(\A)^2}
                       \end{pmatrix}\mathsf{V}.
 $$
Since $\mathsf{U}$ and $\mathsf{V}$ are orthogonal matrices, we therefore have: 
\begin{equation}\label{eq:norm of gradient}
\|\nabla f_n^p(\A) \|_{\S_2^n}^2=4p^2\prod_{r=1}^n \big(1-s_r(\A)^2\big)^{2p}\sum_{k=1}^n\frac{s_k(\A)^2}{\left(1-s_k(\A)^2\right)^2}.
\end{equation}

Consequently,  the numerator in the left hand side of~\eqref{eq:ratio exactly} can be evaluated as follows: 
\begin{align}\label{eq:numerator nabla}
\begin{split}
\int_{B_{\S_\infty^n}}\|\nabla f_n^p(\A)\|_{\S_2^n}^2\ud \A&\stackrel{\eqref{eq:weyl}\wedge \eqref{eq:norm of gradient}}{=} 4p^2c_n\int_{[0,1]^n}\left(\sum_{k=1}^n\frac{s_k^2}{\left(1-s_k^2\right)^2}\right)\prod_{r=1}^n \big(1-s_r^2\big)^{2p}\prod_{i=1}^{n-1}\prod_{j=s+1}^n \big|s_i^2 -s_j^2\big| \ud x\\&\ \ \ \ \stackrel{\eqref{eq:J density def}}{=}\frac{c_nZ_n(p)}{2^n} 4p^2\int_{[0,1]^n} \sum_{k=1}^n \frac{1-x_k}{x_k^2}\J_n^p(x)\ud x\\&\ \ \ \ \ \, =\frac{c_nZ_n(p)}{2^n} 4p^2n \bigg(\int_{[0,1]^n}\frac{\J_n^p(x)}{x_n^2}\ud x -\int_{[0,1]^n}\frac{\J_n^p(x)}{x_n}\ud x\bigg),
\end{split}
\end{align}
where the penultimate step of~\eqref{eq:numerator nabla} makes the same change of variable $s=(\sqrt{1-x_1},\ldots,\sqrt{1-x_n})$ as above.  By combining~\eqref{eq:first change of variable} and~\eqref{eq:numerator nabla} we obtain the following expression for the Rayleigh quotient in~\eqref{eq:ratio exactly}:  
\begin{equation}\label{eq:quotient in terms of jacobi}
\frac{\int_{B_{\S_\infty^n}}\|\nabla f_n^p(\A)\|_{\S_2^n}^2\ud \A}{\int_{B_{\S_\infty^n}} |f_n^p(\A)|^2\ud \A}=4p^2\int_{[0,1]^n} \Big(\sum_{s=1}^n \frac{1-x_s}{x_s^2}\Big)\J_n^p(x)\ud x =
4p^2n \bigg(\int_{[0,1]^n}\frac{\J_n^p(x)}{x_n^2}\ud x -\int_{[0,1]^n}\frac{\J_n^p(x)}{x_n}\ud x\bigg).
\end{equation}

To evaluate the right hand side of~\eqref{eq:quotient in terms of jacobi}, consider the following symmetric polynomial $P\in \Q[x_1,\ldots x_n]$:
\begin{equation}\label{eq:def our P}
P(x)\eqdef \sum_{s=1}^n \prod_{t\in [n]\setminus \{s\}} x_t^2+\frac{2}{3}\sum_{s=1}^{n-1}\sum_{t=s+1}^nx_sx_t\prod_{k\in [n]\setminus \{s,t\}} x_k^2,
\end{equation}
where in~\eqref{eq:def our P}, and throughout, we use the common notation $[n]\eqdef \n$.  Thus,
\begin{equation}\label{eq:integrate P}
\frac{1}{Z_n(p)}\int_{[0,1]^n} P(x)\Vn(x)\prod_{s=1}^n\frac{x_s^{2p-2}}{\sqrt{1-x_s}}\ud x\stackrel{\eqref{eq:J density def}}{=}n\int_{[0,1]^n}\frac{\J_n^p(x)}{x_n^2}\ud x+\frac{n(n-1)}{3}\int_{[0,1]^n}\frac{\J_n^p(x)}{x_{n-1}x_n}\ud x. 
\end{equation}
The final integral that appears in~\eqref{eq:integrate P} can be evaluated as follows:
\begin{equation}\label{eq:n-2}
\int_{[0,1]^n}\frac{\J_n^p(x)}{x_{n-1}x_n}\ud x\stackrel{\eqref{eq:J density def}}{=}\frac{1}{Z_n(p)}\int_{[0,1]^n} \Vn(x)\bigg(\prod_{t=1}^{n-2} x_t\bigg) \prod_{s=1}^n\frac{x_s^{2p-1}}{\sqrt{1-x_s}} \ud x=\frac{(4p+n)(4p+n+1)}{4p(4p+1)},
\end{equation}
where the second equality that appears in~\eqref{eq:n-2} is an instantiation of the Aomoto integral formula~\cite{Aom87} (see also e.g.~\cite[Theorem~8.1.2]{AAR99}). To evaluate the integral that appears in the left hand side of~\eqref{eq:integrate P}, observe that   the polynomial $P$ that is defined in~\eqref{eq:def our P} coincides with the Jack polynomial that is denoted in~\cite{aomoto} by $P^{(1/\gamma)}_\lambda$, where $\gamma=1/2$ and $\lambda$ is the vector in $\Z^n$ whose first $n-1$ coordinates equal $2$ and the whose last coordinate equals $0$; this is seen by substituting the aforementioned parameters $\gamma,\lambda$ into equations~(3.2) and~(3.10) of~\cite{aomoto}. The Kadell integral formula~\cite{Kad97}, which appears in~\cite{aomoto} as equation~(3.13) thus gives:  
\begin{equation}\label{eq:use kadell}
\frac{1}{Z_n(p)}\int_{[0,1]^n} P(x)\Vn(x)\prod_{s=1}^n\frac{x_s^{2p-2}}{\sqrt{1-x_s}}\ud x=\frac{n(n+2)(4p+n)(4p+n-2)}{24p(2p-1)}. 
\end{equation}
The first integral that appears in the right hand side of~\eqref{eq:quotient in terms of jacobi} is therefore evaluated as follows: 
\begin{equation}\label{eq:1/x^2}
\int_{[0,1]^n}\frac{\J_n^p(x)}{x_n^2}\ud x\stackrel{\eqref{eq:integrate P}\wedge \eqref{eq:n-2}\wedge \eqref{eq:use kadell}}{=}\frac{(n+2)(4p+n)(4p+n-2)}{24p(2p-1)}-\frac{(n-1)(4p+n)(4p+n+1)}{12p(4p+1)}.
\end{equation}

Now, the second integral that appears in the right hand side of~\eqref{eq:quotient in terms of jacobi} can be evaluated as follows: 
\begin{equation}\label{eq:1/x}
\int_{[0,1]^n}\frac{\J_n^p(x)}{x_n}\ud x\stackrel{\eqref{eq:J density def}}{=}\frac{1}{Z_n(p)}\int_{[0,1]^n} \Vn(x)\bigg(\prod_{t=1}^{n-1} x_t\bigg) \prod_{s=1}^n\frac{x_s^{2p-1}}{\sqrt{1-x_s}} \ud x=\frac{4p+n}{4p},
\end{equation}
where the second step of~\eqref{eq:1/x} is another instantiation of the Aomoto integral formula. The desired identity~\eqref{eq:ratio exactly}  follows by substituting~\eqref{eq:1/x^2} and~\eqref{eq:1/x} into~\eqref{eq:quotient in terms of jacobi}, and then simplifying the resulting expression.   \qed

\section{Proof of Theorem~\ref{thm:alg rad}}

Observe that every $x=(x_1,\ldots,x_n)\in \R^n$ satisfies the following estimates:
\begin{equation}\label{eq:use Cauchy schwarz}
 \|x\|_\bE\le \sum_{s=1}^n |x_s|\cdot\|e_s\|_\bE\stackrel{\eqref{eq:norms are reasonable}}{\le} e^{n^D}\sum_{s=1}^n |x_s|\le \sqrt{n} e^{n^D} \|x\|_{\ell_{\!\!2}^n},
\end{equation}
where the last step of~\eqref{eq:use Cauchy schwarz} is an application of Cauchy--Schwarz, and also:
\begin{equation}\label{eq:use convexity}
 \|x\|_\bE\stackrel{\eqref{eq:1 symmetric}}{=}  \max_{s\in [n]} \frac{\|x\|_\bE+ \|(-x_1,\ldots,-x_{s-1},x_s,-x_{s+1},\ldots,-x_n)\|_{\bE}}{2}\ge \max_{s\in [n]} |x_s|\cdot\|e_s\|_\bE \stackrel{\eqref{eq:norms are reasonable}}{\ge}   \frac{\|x\|_{\ell_{\!\!2}^n}}{\sqrt{n}e^{n^D}},
\end{equation}
where the penultimate step of~\eqref{eq:use convexity} is an application of the triangle inequality for $\|\cdot\|_\bE$. Therefore,
\begin{equation}\label{eq:circumscribed}
\frac{1}{\sqrt{n} e^{n^D}} B_{\!\ell_{\!\!2}^n}  \stackrel{\eqref{eq:use Cauchy schwarz}}{\subset} B_\bE \stackrel{\eqref{eq:use convexity}}\subset \sqrt{n} e^{n^D}\!\! B_{\!\ell_{\!\!2}^n}. 
\end{equation} 

Using the terminology of~\cite{GLS93}, the inclusions~\eqref{eq:circumscribed} mean that $B_{\bE}$ is well--bounded and centered. Since in Theorem~\ref{thm:alg rad} we are assuming that $B_\bE$ is given by a weak membership oracle, the Judin--Nemirovski\u{i} theorem~\cite{JN76} (see~\cite[Section~4.3]{GLS93}), which is a striking  application of the Shallow-Cut Ellipsoid Method~\cite{JN76} (see~\cite[Section~3.3]{GLS93}), implies it is possible to optimize linear functionals efficiently on $B_\bE$. In particular, since the dual norm $\|\cdot\|_{\bE^*}$ is given by $\|x\|_{\bE^*}=\max_{y\in B_\bE} \langle x,y\rangle$ for every $x\in \R^n$, where $\langle\cdot,\cdot\rangle:\R^n\times \R^n\to \R$ denotes the standard scalar product on $\R^n$, it follows  that   there is  an algorithm  $\mathrm{Eval}_{\bE^*}$ which takes as input $x\in \R^n$ and outputs in oracle time that is polynomial in both $\log(1/\d)$ and $n^D$ a number $\mathrm{Eval}_{\bE^*}\!(x)\in [0,\infty)$ that is guaranteed to satisfy the following estimates:
\begin{equation}\label{eq:alg for E star}
\Big(1-\frac{\d}{2}\Big)\|x\|_{\bE^*}\le \mathrm{Eval}_{\bE^*}\!(x)\le \Big(1+\frac{\d}{2}\Big)\|x\|_{\bE^*}. 
\end{equation}

Define $q=p/(p-1)$ if $p>1$, and $q=\infty$ if $p=1$. By~\cite[Lemma~2]{Gow89}, if we denote
\begin{equation}\label{eq:quadrant net}
K_q^n\eqdef \big\{a=(a_1,\ldots,a_n)\in \R^n:\ \|a\|_{\!\ell_{\!\! q}^n}\le 1\ \mathrm{and}\ a_1\ge a_2\ge \cdots \ge a_n\ge 0\big\},
\end{equation}
then for every $0<\d<1$ there exists a subset $\GG(n,q,\d)$ of $K_{q}^n$ such that:
\begin{equation}\label{gowers net}
\forall x\in K_q^n, \quad d_{\ell_{\!\!q}^n}\big(x,\GG(n,q,\d)\big)= \min_{a\in \GG(n,q,\d)}\|x-a\|_{\!\ell_{\!\!q}^n}\le \frac{\d}{2},
\end{equation}
i.e., $\GG(n,q,\d)$ is $(\d/2)$-dense in $K_q^n$ with respect to the $\ell_{\!\!q}^n$ norm, yet its size satisfies:
\begin{equation}\label{eq:net size}
|\GG(n,q,\d)|\le n^{\frac{2\log\left(\frac{30}{\d}\right)}{\log\left(1+\frac{\d}{6}\right)}}.
\end{equation}
Therefore, by~\eqref{eq:net size} for fixed $\d$ the size of $\GG(n,q,\d)$ grows like a fixed power of $n$, in marked contrast to the  fact (which is a straightforward consequence of volume comparison)  that  every $(\d/2)$-dense subset of a unit ball of an $n$-dimensional normed space  must have size at least $(2/\d)^n$. 

The algorithm $\mathrm{Alg}_p$ of Theorem~\ref{thm:alg rad}  evaluates $\mathrm{Eval}_{\bE^*}\!(a)$ for each $a\in \GG(n,q,\d)$ and outputs the number: 
\begin{equation}\label{eq:def our alg}
\mathrm{Alg}_p(n,D,\d,\bE)\eqdef 2\max_{a\in \GG(n,q,\d)} \mathrm{Eval}_{\bE^*}\!(a). 
\end{equation}
The oracle running time of this simple procedure is a polynomial in both   $\log(1/\d)$ and $n^D$ times the size of $\GG(n,q,\d)$, so by~\eqref{eq:net size} it is at most $n^d$ where the degree $d$ satisfies the desired bound~\eqref{eq:deg}.  Furthermore,
\begin{equation}\label{eq:duality}
\diam_{\ell_{\!\!p}^n}(B_\bE) =2\sup_{x\in B_\bE} \|x\|_{\!\ell_{\!\!p}^n}=2\sup_{x\in B_\bE} \sup_{y\in B_{\!\ell_{\!\!q}^n}}\langle x,y\rangle = 2\sup_{y\in B_{\!\ell_{\!\!q}^n}} \sup_{x\in B_\bE}\langle x,y\rangle  =2\sup_{y\in B_{\!\ell_{\!\!q}^n}}\|x\|_{\bE^*}.
\end{equation}
Since $\GG(n,q,\d)$ is a subset of $B_{\!\ell_{\!\!q}^n}$, this implies in particular that the following a priori upper bound on the output of $\mathrm{Alg}_p$ holds, which gives the second inequality in the desired conclusion~\eqref{eq:guarantees of alg} of Theorem~\ref{thm:alg rad}:
$$
\mathrm{Alg}_p(n,D,\d,\bE)\stackrel{\eqref{eq:alg for E star}\wedge \eqref{eq:def our alg}}{\le} 2\max_{a\in \GG(n,q,\d)} \Big(1+\frac{\d}{2}\Big)\|a\|_{\bE^*}\le 2\Big(1+\frac{\d}{2}\Big) \sup_{y\in B_{\!\ell_{\!\!q}^n}}\|x\|_{\bE^*} \stackrel{\eqref{eq:duality}}{\le} \Big(1+\frac{\d}{2}\Big)\diam_{\ell_{\!\!p}^n}(B_\bE).
$$

Define $\DD(n,q,\d)\subset B_{\!\ell_{\!\! q}^n}$ by:  
\begin{equation}\label{eq:def DD}
\DD(n,q,\d)\eqdef \Big\{\big(\e_1 a_{\pi(1)},\e_2a_{\pi(2)},\ldots,\e_n a_{\pi(n)}\big):\ (\pi,(\e_1,\ldots,\e_n),a)\in S_n\times \{-1,1\}^n\times \GG(n,q,\d)\Big\}.
\end{equation}
Then, $\DD(n,q,\d)$ is $(\d/2)$-dense in  $B_{\!\ell_{\!\!q}^n}$. Indeed,  for each  $x=(x_1,\ldots,x_n)\in B_{\!\ell_{\!\!q}^n}$ fix a permutation $\pi=\pi^x\in S_n$ such that $|x_{\pi(1)}|\ge |x_{\pi(2)}|\ge \ldots \ge |x_{\pi(n)}|$. Denoting $\e_i=\e_i^x=\sign(x_{i})\in \{-1,1\}$ for each $i\in [n]$, the vector $y=(\e_{\pi(1)}x_{\pi(1)},\ldots,\e_{\pi(n)}x_{\pi(n)})=(|x_{\pi(1)}|,\ldots,|x_{\pi(n)}|)$ belongs $K_{q}^n$, by~\eqref{gowers net}. Now, use~\eqref{gowers net} to get $a\in  \GG(n,q,\d)$ with $\|y-a\|_{\!\ell_{\!\!q^n}}\le \d/2$, so that the point $b=(\e_1a_{\pi^{-1}(1)},\ldots,\e_na_{\pi^{-1}(n)})$ belongs to $\DD(n,q,\d)$ by~\eqref{eq:def DD}, and satisfies $\|x-b\|_{\!\ell_{\!\!q}^n}=\|y-a\|_{\!\ell_{\!\!q}^n}\le \d/2$.   

Observe that the convex hull $\conv(\DD(n,q,\d))$ of $\DD(n,q,\d)$ satisfies the following inclusion:
\begin{equation}\label{eq:conve inclusion}
\conv\left(\DD(n,q,\d)\right)\supseteq \Big(1-\frac{\d}{2}\Big)B_{\!\ell_{\!\!q}^n}.
\end{equation}
 In fact, more generally, for any normed space $\X=(\R^n,\|\cdot\|_\X)$, if $\sub\subset B_\X$ is $(\d/2)$-dense in $B_\X$  with respect to the metric induced by $\|\cdot\|_\X$, then $\conv(\sub)\supseteq (1-\d/2)B_\X$.  Indeed, suppose contrapositively that there exists  $x\in (1-\d/2)B_\X\setminus \conv(\sub)$. Then, by the separation theorem we can fix  $y^*\in \R^n$ with $\|y^*\|_{\X^*}=1$ such that $\langle y^*,x\rangle > \langle y^*,c\rangle$ for every $c\in \sub$. By the Hahn--Banach theorem, we can fix $y\in B_\X$ such that $\langle y^*,y\rangle =1$. Finally, because $\sub$ is $(\d/2)$-dense in $B_\X$, we can fix $c_y\in \sub$ such that $\|y-c_y\|_\X\le \d/2$. We now arrive at the desired contradiction as follows: 
\begin{equation}\label{eq:use separating functional}
1-\frac{\d}{2}\ge \|x\|_\X\ge \langle y^*,x\rangle > \langle y^*,c_y\rangle= \langle y^*y\rangle -\langle y^*,c_y-y\rangle \ge 1-\|c_y-y\|_\X\ge 1-\frac{\d}{2},
\end{equation}
 where the second and the penultimate steps of~\eqref{eq:use separating functional} use the fact that $y^*$ has norm $1$ as a linear functional in $\X^*$, the first step of~\eqref{eq:use separating functional} uses the assumption $x\in (1-\d/2)B_\X$, the third step of~\eqref{eq:use separating functional} uses the assumed property of the separating functional $y^*$ for $c=c_y\in \sub$, and the final step of~\eqref{eq:use separating functional} uses the choice of $c_y$.  
  
We can now justify the first inequality of the desired conclusion~\eqref{eq:guarantees of alg} of Theorem~\ref{thm:alg rad} as follows: 
\begin{align}
\nonumber \mathrm{Alg}_p(n,D,\d,\bE) &\stackrel{\eqref{eq:alg for E star}\wedge \eqref{eq:def our alg}}{\ge} 2\max_{a\in \GG(n,q,\d)} \Big(1-\frac{\d}{2}\Big)\|a\|_{\bE^*}\\ & \ \ \ \ \label{eq:finish alg} \stackrel{\eqref{eq:def DD}}{=}2\max_{a\in \DD(n,q,\d)} \Big(1-\frac{\d}{2}\Big)\|a\|_{\bE^*}\\   &  \ \ \ \  \ = 2\Big(1-\frac{\d}{2}\Big)\max_{x\in \conv\left(\DD(n,q,\d)\right)} \|x\|_{\bE^*}\label{eq:use convexity of norm} \\  \nonumber & \ \ \ \  \stackrel{\eqref{eq:conve inclusion}}{\ge} 2\Big(1-\frac{\d}{2}\Big)\max_{x\in \left(1-\frac{\d}{2}\right)B_{\!\ell_{\!\!q}^n}} \|x\|_{\bE^*}\\ \nonumber &\ \ \ \  \  =2\Big(1-\frac{\d}{2}\Big)^2\max_{x\in B_{\!\ell_{\!\!q}^n}} \|x\|_{\bE^*}\\ \nonumber & \ \ \ \  \stackrel{\eqref{eq:duality}}{=}\Big(1-\frac{\d}{2}\Big)^2\diam_{\ell_{\!\!p}^n}(B_\bE)\\ \nonumber &\ \ \ \  \ >(1-\d)\diam_{\ell_{\!\!p}^n}(B_\bE),
\end{align}
where~\eqref{eq:finish alg} is the step in which the assumption that $\|\cdot\|_\bE$ is a symmetric norm is used, as this implies (by definition) that its dual $\|\cdot\|_{\bE^*}$ is also a symmetric norm, and~\eqref{eq:use convexity of norm}  holds since $\|\cdot\|_{\bE^*}$ is a convex function.  \qed

\begin{remark}\label{rem:approx ratios} {\em In the Introduction we stated the contrast between Theorem~\ref{thm:alg rad} and the known nonexistence of any oracle polynomial time algorithm that computes a $O(1)$-factor approximation to the $\ell_p^n$-diameter of an arbitrary (so, not necessarily symmetric, as in Theorem~\ref{thm:alg rad}) convex body in $\R^n$ that is given by a weak membership oracle, where $1\le p<\infty$ is fixed. This contrast is even more severe because the known impossibility results rule out approximation factors that tend to $\infty$ with $n$, as we will next recall.

For  $n\in \N$, $h\ge 3$ and $1\le p<\infty$, let $D_p^h(n)$ be the smallest $\alpha\ge 1$ for which there is a deterministic algorithm that takes as input a convex body $K\subset \R^n$  given by a weak membership oracle, and outputs in oracle time at most $n^h$ a number $\Delta_p(K)$ satisfying $\diam_{\ell_{\!\!p}^n}(K)\le \Delta_p(K)\le \alpha \diam_{\ell_{\!\!p}^n}(K)$. Similarly, let $R_p^h(n)$ denote the smallest $\alpha\ge 1$ for which there is a randomised algorithm that  outputs in oracle time at most $n^h$ a number $\Delta_p(K)$ satisfying $\diam_{\ell_{\!\!p}^n}(K)\le \Delta_p(K)\le \alpha \diam_{\ell_{\!\!p}^n}(K)$ with probability at least, say, $1/3$. 

Using the above  notations, for  $2\le p<\infty$ we have  $D_p^h(n)\asymp_h R_p^h(n)\asymp_h (n/\log n)^{1/p}$; when $p=2$ these bounds on $D_2^h(n)$ are due to~\cite{BF87}, and these bounds for both $D_p^h(n)$ and $R_p^h(n)$, and arbitrary $2<p<\infty$, are due to~\cite{BGKKLM98,BGKKLS01}. If $1\le p<2$, then  $R_p^h(n)\asymp_h \sqrt{n/\log n}$ by~\cite{KN07,KN08}, improving by unbounded polylog factor over the corresponding bounds of~\cite{BGKKLM98,BGKKLS01}. For deterministic algorithms, when $1\le p<2$ we have   $D_p^h(n)\lesssim_h \sqrt{n}/(\log n)^{1-1/p}$ by~\cite{BGKKLM98,BGKKLS01}, and the best-known lower bound is $D_p^h(n)\gtrsim_h \sqrt{n/\log n}$, which follows from the trivial lower bound $D_p^h(n)\ge R_p^h(n)$ and the aforementioned asymptotic evaluation of $R_p^h(n)$. Thus, for fixed $1\le p<2$ there remains a gap of $(\log n)^{1/p-1/2}$, which tends to $\infty$ with $n$, between the best-known upper and lower bounds on the deterministic approximation ratio $D_p^h(n)$. While it seems reasonable to expect that this gap (noted in~\cite{BGKKLM98,BGKKLS01,Bri02}) could be closed using available methods, to the best of our knowledge this  has not yet been carried out, even though it would be worthwhile to do so. 
}
\end{remark}

\section{Entires in the complex numbers and the quaternions}\label{section:other entries}

Here we will show how to deduce Theorem~\ref{thm:weak revers iso unitarily} for unitarily invariant norms on  $\M_n(\mathbb{C})$ and $\M_n(\mathbb{H})$ from its statement about such norms on $\M_n(\R)$ that we already proved. This will  imply Theorem~\ref{thm:sep formula} and Theorem~\ref{thm:spectral} for complex and quaternionic entries by the reductions from~\cite{naor2022} that were recalled in the Introduction. Suppose that $\bE=(\R^n,\|\cdot\|_\bE)$ is a symmetric norm on $\R^n$. To not supress the underlying entries, we will henceforth denote by $\S_\bE^n(\R), \S_\bE^n(\C), \S_\bE^n(\HH)$ the corresponding unitary ideals on $\M_n(\R), \M_n(\C),\M_n(\HH)$, respectively. Thus, $\S_\bE^n(\R)$ coincides with the space that  we previously  denoted $\S_\bE^n$, and the dimensions over $\R$ of $\S_\bE^n(\R), \S_\bE^n(\C), \S_\bE^n(\HH)$ satisfy $\dim_\R(\S_\bE^n(\R))=n^2,\dim_\R(\S_\bE^n(\C))=2n^2, \dim_\R(\S_\bE^n(\HH))=4n^2$.

For complex entries, our goal is to deduce from Theorem~\ref{thm:weak revers iso unitarily}  that there exists an origin-symmetric convex body $L_\C\subset \M_n(\C)$ that satisfies the following three properties: 
\begin{equation}\label{eq:weak revese iso for unitarily complex}
L_\C\subset B_{\S_\bE^n(\C)}\qquad \mathrm{and}\qquad \vol_{2n^2}(L_\C)^{\frac{1}{2n^2}}\gtrsim \vol_{2n^2}(B_{\S_\bE^n(\C)})^{\frac{1}{2n^2}}\qquad\mathrm{and }\qquad  \iq(L_\C)\lesssim n. 
\end{equation} 
We will henceforth use the natural identification 
\begin{equation}\label{eq:indentitification}
\M_n(\C)\cong \M_n(\R)\times \M_n(\R),
\end{equation}
 which is obtained by assigning to each complex matric $\A\in \M_n(\C)$ the pair $(\Re(\A),\Im(\A))\in \M_n(\R)\times \M_n(\R)$ of real matrices, where $\Re(\A),\Im(\A)$ denote the entry-wise real part and imaginary part of $\A$, respectively. 
 
Every $\A\in \M_n(\C)$ satisfies: 
\begin{multline*}
 \|\A\|_{\S_\bE^n(\C)} =\|\Re(\A)+i\Im(\A)\|_{\S_\bE^n(\C)}  \le \|\Re(\A)\|_{\S_\bE^n(\C)} +\|\Im(\A)\|_{\S_\bE^n(\C)}\\ = \|\Re(\A)\|_{\S_\bE^n(\R)} +\|\Im(\A)\|_{\S_\bE^n(\R)}\le 2\max\big\{ \|\Re(\A)\|_{\S_\bE^n(\R)},\|\Im(\A)\|_{\S_\bE^n(\R)} \big\}.
\end{multline*}
 Equivalently, under the identification~\eqref{eq:indentitification} the following inclusion holds:  
\begin{equation}\label{eq:first inclusion}
B_{\S_\bE^n(\R)} \times B_{\S_\bE^n(\R)}\subset 2B_{\S_\bE^n(\C)}. 
\end{equation}
 Conversely, letting  $\overline{\A}$ denote the entry-wise complex conjugate of $\A\in \M_n(\C)$, we have
$$
\forall \A\in \M_n(\C),\qquad \|\Re(\A)\|_{\S_\bE^n(\R)}=\frac{1}{2} \|{\A+\overline{\A}}\|_{\mathsf{S}_{\mathbf{E}}^n(\mathbb{R})}=\frac{1}{2} \|{\A+\overline{\A}}\|_{\mathsf{S}_{\mathbf{E}}^n(\mathbb{C})} \leq \frac{\|\A\|_{\mathsf{S}_{\mathbf{E}}^n(\mathbb{C})} +\|\overline{\A}\|_{\mathsf{S}_{\mathbf{E}}^n(\mathbb{C})} }{2}= \|\A\|_{\mathsf{S}_{\mathbf{E}}^n(\mathbb{C})},
$$
where  the last step uses the fact that $\A$ and $\overline{\A}$ have the same singular values, whence 
 $\|\A\|_{\mathsf{S}_{\mathbf{E}}^n(\mathbb{C})}= \|\overline{\A}\|_{\mathsf{S}_{\mathbf{E}}^n(\mathbb{C})}$. For the same reason, we also have the following estimate: 
$$
\forall \A\in \M_n(\C),\qquad \|\Im(\A)\|_{\S_\bE^n(\R)}=\frac{1}{2} \|{-i\A+i\overline{\A}}\|_{\mathsf{S}_{\mathbf{E}}^n(\mathbb{R})}=\frac{1}{2} \|{-i\A+i\overline{\A}}\|_{\mathsf{S}_{\mathbf{E}}^n(\mathbb{C})} \leq \frac{\|\A\|_{\mathsf{S}_{\mathbf{E}}^n(\mathbb{C})} +\|\overline{\A}\|_{\mathsf{S}_{\mathbf{E}}^n(\mathbb{C})} }{2}= \|\A\|_{\mathsf{S}_{\mathbf{E}}^n(\mathbb{C})}.
$$
We have thus checked that:
$$
\forall \A\in \M_n(\C),\qquad \max \big\{\|\Re(\A)\|_{\S_\bE^n(\R)},\|\Im(\A)\|_{\S_\bE^n(\R)}\big\}\le \|\A\|_{\mathsf{S}_{\mathbf{E}}^n(\mathbb{C})}. 
$$
Equivalently, under the identification~\eqref{eq:indentitification} the following inclusion holds:  
\begin{equation}\label{eq:second inclusion}
B_{\S_\bE^n(\C)}\subset B_{\S_\bE^n(\R)} \times B_{\S_\bE^n(\R)}. 
\end{equation}

Now, an application of Theorem~\ref{thm:weak revers iso unitarily} yields an origin-symmetric convex body $L_\R\subset \M_n(\R)$ that satisfies: 
\begin{equation}\label{eq:weak revese iso for unitarily real}
L_\R\subset B_{\S_\bE^n(\R)}\qquad \mathrm{and}\qquad \vol_{n^2}(L_\R)^{\frac{1}{n^2}}\gtrsim \vol_{n^2}(B_{\S_\bE^n(\R)})^{\frac{1}{n^2}}\qquad\mathrm{and}\qquad  \iq(L_\R)\lesssim n. 
\end{equation} 
Under the identification~\eqref{eq:indentitification}, define $L_\C\subset \M_n(\C)$ by: 
\begin{equation}\label{eq:def LC}
L_\C\eqdef \frac12 \big(L_\R\times L_\R\big)\stackrel{\eqref{eq:weak revese iso for unitarily real}}{\subset}\frac12 \big(B_{\S_\bE^n(\R)}\times B_{\S_\bE^n(\R)} \big)\stackrel{\eqref{eq:first inclusion}}{\subset} B_{\S_\bE^n(\C)}.  
\end{equation}
Thus, the first part of~\eqref{eq:weak revese iso for unitarily complex} holds. The second part of~\eqref{eq:weak revese iso for unitarily complex} also holds because: 
$$
\vol_{2n^2}(L_\C)^{\frac{1}{2n^2}}\stackrel{\eqref{eq:def LC}}{=}\frac12 \vol_{n^2}(L_\R)^{\frac{1}{n^2}}\stackrel{\eqref{eq:weak revese iso for unitarily real}}{\gtrsim} \vol_{n^2}(B_{\S_\bE^n(\R)})^{\frac{1}{n^2}}=\vol_{2n^2} (B_{\S_\bE^n(\R)}\times B_{\S_\bE^n(\R)})^{\frac{1}{2n^2}}\stackrel{\eqref{eq:second inclusion}}{\ge}\vol_{2n^2} (B_{\S_\bE^n(\C)})^{\frac{1}{2n^2}}.
$$
Finally, the third part of~\eqref{eq:weak revese iso for unitarily complex}  is justified as follows, thus completing the treatment of complex entries: 
\begin{equation*}
    \iq(L_\C)=\stackrel{\eqref{eq:def LC}}{=}\frac{\vol_{2n^2-1}\big(\partial(L_\R\times L_\R)\big)}{\vol_{2n^2}(L_\R\times L_\R)^{1-\frac{1}{2n^2}}}
    = \frac{\vol_{2n^2-1}\big((\partial L_\R) \times L_\R\big)+\vol_{2n^2-1}\big(L_\R\times (\partial L_\R) \big)}{\vol_{n^2}(L_\R)^{2-\frac{1}{n^2}}}=2\iq(L_\R) \stackrel{\eqref{eq:weak revese iso for unitarily real}}{\lesssim} n. 
\end{equation*}

In the case of quaternionic entries one applies mutatis mutandis the analogous reasoning by writing uniquely each $\A\in \M_n(\HH)$ as $\A=\mathsf{B}+i\mathsf{C}+j\mathsf{D}+k\mathsf{E}$ for $\mathsf{B},\mathsf{C},\mathsf{D},\mathsf{E}\in \M_n(\R)$ that satisfy the following bounds:
\begin{equation*}
    \max\left\{\|\mathsf{B}\|_{\mathsf{S}_\mathbf{E}^n(\mathbb{R})}, \|\mathsf{C}\|_{\mathsf{S}_\mathbf{E}^n(\mathbb{R})}, \|\mathsf{D}\|_{\mathsf{S}_\mathbf{E}^n(\mathbb{R})}, \|\mathsf{E}\|_{\mathsf{S}_\mathbf{E}^n(\mathbb{R})}\right\}\leq \|\A\|_{\mathsf{S}_\mathbf{E}^n(\mathbb{H})}\leq 4   \max\left\{\|\mathsf{B}\|_{\mathsf{S}_\mathbf{E}^n(\mathbb{R})}, \|\mathsf{C}\|_{\mathsf{S}_\mathbf{E}^n(\mathbb{R})}, \|\mathsf{D}\|_{\mathsf{S}_\mathbf{E}^n(\mathbb{R})}, \|\mathsf{E}\|_{\mathsf{S}_\mathbf{E}^n(\mathbb{R})}\right\}.
\end{equation*}

\bigskip

\subsection*{Acknowledgements} We are grateful to Alexei Borodin for helpful discussions, and to Ofer Zeitouni for pointers to the literature.

\bibliographystyle{abbrv}

\bibliography{SEP-ideals}

\end{document}